\documentclass[12pt, twoside]{article}
\usepackage{amsmath,amsthm}
\usepackage{amsfonts}
\usepackage{amssymb,latexsym}
\usepackage{enumerate}
\usepackage{url}
\usepackage[colorlinks=true,linkcolor=blue,urlcolor=blue]{hyperref}
\newtheorem{theorem}{Theorem}[section]
\newtheorem{remark}[theorem]{Remark}
\newtheorem{proposition}[theorem]{Proposition}
\newtheorem{lemma}[theorem]{Lemma}
\newtheorem{definition}[theorem]{Definition}

\def\ind{1{\hskip -3 pt}\hbox{\textsc{I}}}
\def\n{\noindent}
\def\n{\noindent}

\def\Om{\Omega}
\def\E{\mathcal E}

\def\F{\mathcal F}

\def\n{\noindent}

\def\E{\mathcal E}

\def\C{\mathbb{C}^n}

\numberwithin{equation}{section}
\begin{document}
\setlength{\baselineskip}{18truept}
\pagestyle{myheadings}

\markboth{Nguyen Van Phu }{Subsolution theorem in weighted energy classes}
\title { A note on the subsolution theorem in weighted energy classes of $m$-subharmonic functions  with given boundary values}
\author{
	Nguyen Van Phu* 
	\\ *Faculty of Natural Sciences, Electric Power University,\\ Hanoi,Vietnam.\\
	\\E-mail: phunv@epu.edu.vn\\ Orcid  \url{https://orcid.org/0000-0002-2851-7250}}

\date{}
\maketitle

\renewcommand{\thefootnote}{}



\renewcommand{\thefootnote}{\arabic{footnote}}
\setcounter{footnote}{0}

\begin{abstract}
 In this note, we deal with the existence of solutions of the weighted complex $m$-Hessian equation $-\chi(u)H_{m}(u)=\mu$ in the class $\mathcal{E}_{m,\chi}(f,\Omega)$  if there exists a subsolution in this class, where the given boundary value $f\in\mathcal{E}_m(\Omega)\cap MSH_m(\Omega).$ This is a generalization of the result in the paper \cite{PDtaiwan} where we proved that the subsolution theorem is true in the class $\mathcal{E}_{m,\chi}(f,\Omega)$ in the case when the given boundary value $f\in\mathcal{N}_m(\Omega)\cap MSH_m(\Omega).$
\end{abstract}

\section{Introduction and the main result}
The complex Monge-Amp\`ere operator $(dd^c.)^n$ plays a central role in pluripotential theory and has attracted the interest of mathematicians in recent years. In 1982, in their famous and foundational mathematical work, E. Bedford and B. A. Taylor \cite{BT1} demonstrated that this operator  is well-defined in the class of locally bounded plurisubharmonic functions over a bounded domain $\Omega$ in  $\mathbb{C}^n$ and its image lies in the class of non-negative measures. After that, in his seminal work, U. Cegrell \cite{Ce98, Ce04} introduced several classes of functions  $\mathcal{F}(\Omega)$ and $\mathcal{E}(\Omega)$ which are not necessarily locally bounded, on which the complex Monge-Ampère operator is well-defined. In addition to investigating the existence of the Monge-Ampère operator, solving the Monge-Amp\`ere equation  $(dd^c.)^n=\mu$ in the class of plurisubharmonic functions constitutes a central problem in the pluripotential theory - this is  the tool for constructing metrics with preassigned Ricci curvature. In 2004, U. Cegrell \cite{Ce04} established a solution to the aforementioned equation under the assumption that the measure  $\mu$ vanishes on pluripolar sets. Subsequently, in 2010, R. Czyż \cite{Cz10} addressed the weighted Monge-Amp\`ere equation   $-\chi(u)(dd^cu)^n=\mu$ under the same condition of measure   $\mu$.
Following that, L. M.  Hai and P. H. Hiep \cite{HH11} further generalized Czyż’s result to a broader class of Monge-Ampère type equations   $\chi(u(z),z)(dd^c.)^n=\mu.$ In 1995,  S. Ko{\l}odziej  \cite{Ko95} established the existence of bounded solutions to the Monge-Amp\`ere equation under the assumption that it has a bounded subsolution. This result is widely known as Kołodziej’s subsolution theorem. Specifically, if $\mu\leq (dd^cw)^n,$ where $w\in PSH^-(\Omega)\cap L^{\infty}(\Omega),$ then there exists $u\in PSH^-(\Omega)\cap L^{\infty}(\Omega)$ such that $\mu=(dd^cu)^n.$ In 2009, P. Ahag, U. Cegrell, R. Czyż, and P. H. Hiep \cite{ACCH} proved the subsolution theorem for the class of unbounded plurisubharmonic functions $\mathcal{E}(\Omega).$ In the sequel,  L. M. Hai, P. H. Hiep, N. X. Hong and N. V Phu \cite[Theorem 3.1]{HHHP} proved the subsolution theorem for the weighted complex Monge-Amp\`epe equation $\mu=-\chi(u)(dd^cu)^n$ in the class $\mathcal{E}_{\chi}(\Omega).$  Note that, in  \cite{Ko95}, the measure $\mu$ vanishes on pluripolar sets, whereas in \cite{ACCH} and \cite{HHHP}, the condition that measure $\mu$ vanishes on all pluripolar sets of $\Omega$ is not required; it only needs to have finite mass. Let $f\in\mathcal{E}(\Omega)$ be a maximal plurisubharmonic function and measure $\mu$ vanishes on all pluripolar sets of $\Omega$, Benelkourchi \cite{Be13} proved that the complex Monge-Amp\`ere type equation $F(u,z)\mu=(dd^cu)^n$ has solution in $\mathcal{N}^a(f,\Omega)$ if it has subsolution in $\mathcal{N}^a(\Omega).$ Namely, if we have  $F(w,z)\mu\leq (dd^cw)^n$ for some $w\in \mathcal{N}^a(\Omega)$  then there exists a function $u$ such that $F(u,z)\mu=(dd^cu)^n$ and  $f+w\leq u\leq f.$  \\
B{\l}ocki \cite{Bl1} and  Sadullaev, Abdullaev \cite{SA12} introduced $m$-subharmonic functions which are extensions of the plurisubharmonic functions. They also proved that the complex $m$-Hessian operator $H_m(.) = (dd^c.)^m\wedge \beta^{n-m}$ is well-defined in the class of locally bounded $m$-subharmonic functions. L. H. Chinh \cite{Ch12} introduced the Cegrell classes $\mathcal{F}_m(\Omega)$ and $\mathcal{E}_m(\Omega)$ which are not necessarily locally bounded and  the complex $m$-Hessian operator  is well defined in these classes. In the theory of $m$-subharmonic functions, the complex $m$-Hessian equation $\mu=H_m(u)$ plays an important role. Mathematicians are interested in solving the $m$-Hessian equation when it has a subsolution. In \cite{Cu13}, N. N. Cuong proved the subsolution theorem for the $m$-Hessian equation in the class of bounded $m$-subharmonic functions. After that, V. V. Hung and N. V. Phu \cite{HP17} proved that the subsolution theorem is true in the class $\mathcal{E}_{m}(\Omega).$  In 2023, D. T. Duong and N. V. Thien \cite{DT23} introduced and investigated the weighted energy class of $m$-subharmonic functions $\mathcal{E}_{m,\chi}(\Omega).$  Following that,  N. V. Phu and N. Q. Dieu \cite{PD2023} proved the subsolution theorem for the weighted complex $m$-Hessian equation $-\chi(u)H_m(u)=\mu$ in the class $\mathcal{E}_{m,\chi}(\Omega).$ Note that, in \cite{Cu13}, the measure $\mu$ puts no mass on all $m$-polar sets, whereas in \cite{HP17} and \cite{PD2023}, the condition that measure $\mu$ vanishes on all $m$-polar sets of $\Omega$ is not required. In the case when the measure $\mu$ puts no mass on all $m$-polar set, L. M. Hai and V. V. Quan \cite[Theorem 1.1]{HQ21} solved the complex $m$-Hessian type equation $H_m(u,z)=F(u,z)d\mu$ in the class $\mathcal{D}_m(\Omega)$ if it has a subsolution in this class, where the function $t\mapsto F(t,z)$ is continuous and non-decreasing. After that, N. V. Phu and N. Q. Dieu \cite[Theorem 4.1]{PDjmaa} solved the equation $\chi(u,z)H_m(u)=\mu$ in the class $\mathcal{F}_m^a(f,\Omega)$ without requiring that the function $\chi(t,z)$ is monotone in its first variable for all $z\in\Omega.$\\
\n Let $f\in\mathcal{E}_m(\Omega)$ be a $m$-subharmonic function. Extending the result of \cite{HP17}, Gasmi \cite{Gasmi} proved that if the complex $m$-Hessian equation has a subsolution in the class $\mathcal{N}_{m}(\Omega)$ then it has a solution in the class $\mathcal{N}_{m}(f)$.    Based on the work of Benelkourchi \cite{Be13}, Amal, Asserda and Gasmi \cite{AAG20}   solved the complex $m$-Hessian type equation $H_{m}(u)=F(u,z)d\mu$ in the class $\mathcal{N}_{m}(f)$ if there exists a subsolution in the class $\mathcal{N}_{m}(\Omega).$ In contrast to Benelkourchi’s result for plurisubharmonic functions, the result of Amal, Asserda and Gasmi \cite{AAG20}  is obtained without imposing the condition that the measure $\mu$ vanishes on $m$-polar sets. Developing further the result of Amal, Asserda and Gasmi \cite{AAG20}, N. V. Phu and N. Q. Dieu \cite{PDtaiwan} proved the subsolution theorem for the weighted Hessian equation $-\chi(u)H_m(u)=\mu$ in the class $\mathcal{E}_{m,\chi}(f,\Omega)$ in the case when $f\in\mathcal{N}_m(\Omega)\cap MSH_m(\Omega).$ Note that, similar to Benelkourchi’s result for plurisubharmonic functions, the result of Amal, Asserda and Gasmi \cite{AAG20} also provides an estimate relating   a solution $u$ and a subsolution $w$ via the inequality $f+w\leq u\leq f.$ The result of Phu–Dieu \cite[Theorem 3.3]{PDtaiwan} gives a sharper estimate, namely $w\leq u\leq f$. However, the result of N. V. Phu and N. Q. Dieu \cite{PDtaiwan} is obtained only in the case where the function $f$ belongs to the class $\mathcal{N}_m(\Omega),$ whereas the result of Amal, Asserda and Gasmi \cite{AAG20} is achieved for the function $f$ belonging to the class $\mathcal{E}_m(\Omega)$ which is more general than $\mathcal{N}_m(\Omega).$  Naturally, we aim to extend the result in \cite{PDtaiwan} to the function $f$ belonging to the class $\mathcal{E}_m(\Omega).$ To achieve this, we need a completely different approach from the method used in the paper \cite{PDtaiwan}. In our paper, we will use the recent result on the envelope of $m$-subharmonic functions from  \cite{Pmalay} to obtain the desired result. The main result is contained in the following theorem.
   \begin{theorem}
   	\label{main theorem}
   	Let  $\mu$ be a non-negative finite measure on $\Omega$ and $\chi: [-\infty,0]\to [-\infty,0] $ be a  nondecreasing continuous function such that $\chi(t)<0$ for all $t<0$. Assume that there exists a function $w\in\mathcal E_{m,\chi}(f,\Omega)$ with $\mu\leq -\chi(w)H_{m}(w)$, where the given boundary $f\in\mathcal{E}_m(\Omega)\cap MSH_m(\Omega)$. Then there exists a function $u\in\mathcal{E}_{m,\chi}(f,\Omega)$ such that $u\geq w$ and $$-\chi(u)H_{m}(u)=\mu.$$
   \end{theorem}
  \begin{remark}
  	It should be emphasized that the result of  Amal, Asserda and Gasmi   requires the use of the Schauder–Tikhonov fixed point theorem (see \cite[Theorem 3.1]{AAG20}). This method was previously employed by Benelkourchi \cite{Be13} for plurisubharmonic functions. Compared with the result of Amal, Asserda and Gasmi \cite{AAG20}, we obtain a sharper estimate on the relationship between the solution and the subsolution, using a different and more concise proof, without employing the Schauder–Tikhonov fixed point theorem.
  	\end{remark}
{\bf Acknowledgments.}  
The work is supported by Ministry of Education and Training, Vietnam under the Grant number B2025-CTT-10.\\
This work is written during my visit in Vietnam Institute for Advanced Study in Mathematics (VIASM) during the Spring of 2025. I acknowledge greatfully the  financial support and hospitality VIASM.	

\section{Preliminaries}
\n Some elements of the theory of $m$-subharmonic functions and the complex $m$-Hessian operator can be found in \cite{Bl1}, \cite{SA12}, \cite{Ch12}, \cite{Ch15}, \cite{DiKo} and \cite{T19}.
We denote by $\beta= dd^c\|z\|^2$ the canonical K\"ahler form of $\mathbb{C}^n$ with the volume element $dV_{2n}= \frac{1}{n!}\beta^n$ where $d= \partial +\overline{\partial}$ and $d^c =\frac{\partial - \overline{\partial}}{4i}$. \\
\n Let $1\leq m\leq n$ be integers. Let $\Omega$ be a bounded $m$-hyperconvex domain in $\mathbb{C}^n$, which mean there exists a negative continuous $m$-subharmonic function $\rho:\Omega\to (-\infty,0)$ such that  $\{z\in\Omega:\rho(z)<c\}\Subset\Omega$ for every $c<0.$ Such a function $\rho$ is called the exhaustion function on $\Omega.$ Throughout this paper $\Omega$ will denote a bounded $m$-hyperconvex domain in $\mathbb{C}^{n}.$\\
\n We recall the classes $\mathcal{E}^0_m(\Omega)$, $\mathcal{F}_m(\Omega)$ and $\mathcal{E}_m(\Omega)$ introduced in \cite{Ch12}.
\begin{definition}{\rm
		$$
		\mathcal{E}^0_m=\mathcal{E}^0_m(\Omega)=\{u\in SH^{-}_m(\Omega)\cap{L}^{\infty}(\Omega):
		\underset{z\to\partial{\Omega}}
		\lim u(z)=0, \int\limits_{\Omega}(dd^c u)^m\wedge\beta^{n-m} <\infty\},$$
		$$\mathcal{F}_m=\mathcal{F}_m(\Omega)=\big\{u\in SH^{-}_m(\Omega): \exists \mathcal{E}^0_m\ni u_j\searrow u, \underset{j}\sup\int\limits_{\Omega}(dd^c u_j)^m\wedge\beta^{n-m}<\infty\big\},$$
		and
		\begin{align*}
			&\mathcal{E}_m=\mathcal{E}_m(\Omega)=\big\{u\in SH^{-}_m(\Omega):\forall z_0\in\Omega, \exists \text{ a neighborhood } \omega\ni z_0,  \text{ and }\\
			&\hskip4cm \mathcal{E}^0_m\ni u_j\searrow u \text{ on } \omega, \underset{j}\sup\int\limits_{\Omega}(dd^c u_j)^m\wedge\beta^{n-m} <\infty\big\}.
	\end{align*}}
\end{definition}
\noindent From  Theorem 3.14 in \cite{Ch12} it follows that if $u\in \mathcal{E}_m(\Omega)$, the complex $m$-Hessian $H_m(u)= (dd^c u)^m\wedge\beta^{n-m}$ is well defined and it is a Radon measure on $\Omega$. On the other hand, by Remark 3.6 in \cite{Ch12} the following description of  $\mathcal{E}_m(\Omega)$ may be given
\begin{align*}
	&\mathcal{E}_m=\mathcal{E}_m(\Omega)=\big\{u\in SH^{-}_m(\Omega):\forall U\Subset\Omega, \exists {v\in\mathcal{F}_m(\Omega)}, {v=u}  \text{ on }  {U}\bigl\}.
\end{align*}
\begin{definition} {\rm We say that an $m$-subharmonic function $u$ is $m$-maximal if for every relatively compact open set $K$ on $\Omega$ and for each upper semicontinuous function  $v$ on $\overline{K},$ $v\in SH_{m}(K)$ and $v\leq u$ on $\partial K,$ we have $v\leq u$ on $K.$ The family of $m$-maximal $m$-subharmonic function defined on $\Omega$ will be denoted by $MSH_{m}(\Omega).$ 
	}
\end{definition}

\n According to Theorem 1.2 in \cite{Bl1} we see that $u\in\mathcal{E}_m(\Omega)$ is a $m$-maximal if and only if $H_m(u)=0.$\\
\n Next, we recall the definition of the class $\mathcal{N}_m(\Omega)$ as in Subsection 2.6 in \cite{Pjmaa}.
\begin{definition}{\rm Let $u\in SH_{m}(\Omega),$ and let ${\Omega_{j}}$ be a fundamental sequence of $\Omega,$ which means
		$\Om_j$ is strictly $m-$hyperconvex domain, $\Omega_{j}\Subset\Omega_{j+1}$ and $\cup_{j=1}^{\infty}\Omega_{j}=\Omega.$ Set 
		$$u^{j}(z)=\big(\sup\{\varphi(z):\varphi\in SH_{m}(\Omega), \varphi\leq u \,\,\text{on}\,\, \Omega_{j}^{c}\}\big)^{*},$$
		where $\Omega_{j}^{c}$ denotes the complement of $\Omega_{j}$ on $\Omega.$}
\end{definition}

We can see that $u^{j}\in SH_{m}(\Omega)$ and $u^{j}=u$ on $\Omega_{j}^{c}.$ From definition of $u^{j}$ we see that $\{u^{j}\}$ is an increasing sequence and therefore $\lim\limits_{j\to\infty}u^{j}$ exists everywhere except on an $m$-polar subset on $\Omega.$ Hence, the function $\tilde{u}$ defined by $\tilde{u}=\big(\lim\limits_{j\to\infty}u^{j}\big)^{*}$
is an $m$-subharmonic function on $\Omega.$ Obviously, we have $\tilde{u}\geq u.$ Moreover, if $u\in \mathcal{E}_{m}(\Omega)$ then $\tilde{u}\in  \mathcal{E}_{m}(\Omega)$ and $\tilde{u}\in MSH_{m}(\Omega).$
\begin{definition}{\rm
		Set
		$$\mathcal{N}_{m}=\mathcal{N}_{m}(\Omega)=\{u\in\mathcal{E}_{m}(\Omega): \tilde{u} =0\}.$$
	}
\end{definition}

We have the following inclusion
$$\mathcal{F}_{m}(\Omega)\subset \mathcal{N}_{m}(\Omega) \subset \mathcal{E}_{m}(\Omega).$$

\n Let $\mathcal{K}$ be one of the classes $\mathcal{E}_{m}^{0}(\Omega), \mathcal{F}_{m}(\Omega), \mathcal{N}_{m}(\Omega), \mathcal{E}_{m}(\Omega).$ According to Lemma 2.6 in \cite{DT23}, $\mathcal{K}$ is a convex cone, i.e., if $u,v\in \mathcal{K}$ then $au+bv\in \mathcal{K}$ for arbitrary nonnegative constants $a,b.$ Moreover, if $u\in \mathcal{K}, v\in SH_m^-(\Omega)$ then $\max(u,v)\in \mathcal{K}.$ \\
\n  We say that a $m$-subharmonic function defined on $\Omega$ belongs to the class $\mathcal{K}(f,\Omega)$, where $f\in \mathcal{E}_{m}(\Omega)$ if there exists a function $\varphi\in\mathcal{K}$ such that 
$$ f\geq u\geq f+\varphi.$$
Note that $\mathcal{K}(0,\Omega)=\mathcal{K}.$ 

\n  We recall some results on weighted $m$-energy classes in \cite{DT23}.
\begin{definition}{\rm
		Let $\chi: [-\infty,0]\longrightarrow [-\infty,0]$ be an increasing continuous function. We set
		$$\mathcal{E}_{m,\chi}(\Omega)=\{u\in SH_{m}(\Omega):\exists(u_{j})\in\mathcal{E}_{m}^{0}(\Omega), u_{j}\searrow u, \sup_{j}\int_{\Omega}(-\chi)\circ u_{j}H_{m}(u_{j})<+\infty\}.$$
	}
\end{definition}

\n Note that in the case when $\chi\equiv -1,$ then the weighted $m$-energy classes $\mathcal{E}_{m,\chi}(\Omega)$ is the Cegrell energy  classes $\mathcal{F}_{m}(\Omega).$ 
According Theorem 3.3 in \cite{DT23}, if $\chi\not\equiv 0$ then $\mathcal{E}_{m,\chi}(\Omega)\subset \mathcal{E}_{m}(\Omega)$ which means that the complex $m$-Hessian operator is well - defined on class $\mathcal{E}_{m,\chi}(\Omega)$ and if $\chi(t)<0$ for all $t<0$ then we have $\mathcal{E}_{m,\chi}(\Omega)\subset \mathcal{N}_{m}(\Omega).$ \\
\begin{definition}\label{dn61}{\rm Let $E\subset \Omega$ be a Borel subset. The $m$-capacity of $E$ with respect to $\Omega$ is defined in \cite{Ch15} by
		$$C_{m}(E)=C_m(E,\Omega)= \sup\Bigl\{\int\limits_{E}H_m(u): u\in SH_m(\Omega), -1\leq u\leq 0\Bigl\}.$$}
\end{definition}

As in \cite{Ch12}, we have the following Definition.
\begin{definition}{\rm
		A sequence $u_j\in SH_m^-(\Omega)$ converges to $u$ in $C_m$-capacity if
		$$C_m(K\cap\{|u_j-u|>\delta\})\to 0\ \text{as}\ j\to+\infty,\ \forall\ K\subset\subset\Omega,\ \forall\ \delta>0.$$}
\end{definition}
\n According to Lemma 2.9 in \cite{T19}, if a sequence of $m$-subharmonic functions $\{u_j\}$ converges monotonically to a $m$-subharmonic function $u$ then  $u_j\to u$ in $C_m$ as $j\to\infty.$ \\
\n For the convenience of the readers, we recall some results that will be used in this paper.
\begin{proposition}[Corollary 3.3 in  \cite{PD2023}]\label{Pdebrecen}
	Let $\chi:\mathbb R^-\longrightarrow\mathbb R^-$ be an increasing continuous function with $\chi(-\infty) >-\infty$. Let $\{u_j, u\}\subset\mathcal E_{m}(\Omega)$, be such that $u_j\geq v$, $\forall j\geq 1$ for some $v\in\mathcal E_{m}(\Omega)$ and that $u_j \to u \in \mathcal E_{m}(\Omega)$ in $C_m$. Then 
	$-\chi (u_j)H_{m}(u_{j})\to -\chi (u)H_{m}(u)$ weakly.
\end{proposition}	
\n Let $u\in\mathcal{E}_m(\Omega).$ We have $H_m(u)=\ind_{\{u=-\infty\}}H_m(u)+\ind_{\{u>-\infty\}}H_m(u).$
By  Lemma 2.16 in \cite{Gasmi} (see also Theorem 2.8 in \cite{DT23}), we see that the measure $\ind_{\{u>-\infty\}}H_m(u)$ puts no mass on $m$-polar sets. According to Theorem 5.3 in \cite{Ch15}, there exists $\varphi\in\mathcal{E}_m^0(\Omega)$ and $0<f\in L_{loc}^1(H_m(\varphi))$ such that $\ind_{\{u>-\infty\}}H_m(u)=fH_m(\varphi).$ Therefore, we have Cegrell's decomposition theorem $H_m(u)=fH_m(\varphi)+\nu_u,$ where $\nu_u=\ind_{\{u=-\infty\}}H_m(u)$ is a positive measure carried by an $m$-polar set (see also Theorem 2.15 in \cite{Gasmi}). Note that the $m$-polar part of the Hessian measure $H_m(u)$ is supported by $\{u=-\infty\}.$ Using it we have the following Proposition.
\begin{proposition}[Proposition 2.9 in \cite{Pjmaa}]\label{pro2.9}
	\n  Assume that $u,v,u_k\in\E_m(\Om), k=1,\cdots,m-1$ with $u\geq v $ on $\Om$ and  $T=dd^cu_1\wedge\cdots\wedge dd^cu_{m-1}\wedge\beta^{n-m}.$ Then we have
	$$\ind_{\{u=-\infty\}}dd^cu\wedge T\leq \ind_{\{v=-\infty\}}dd^cv\wedge T.$$
	In particular, if $u,v\in\E_m(\Om)$ are such that $u\geq v$ then for every $m$-polar set $A\subset\Om$ we have
	$$ \int_{A}H_m(u)\leq\int_{A}H_m(v).$$
\end{proposition}
\n We recall the following result concerning the $m$-polar part of the Hessian measure. It was proved in \cite[Lemma 4.12]{ACCH} for psh functions.
\begin{proposition}[ Lemma 5.1 in \cite{Gasmi}]\label{pluripolarpart}
	Let $u,v\in\mathcal{E}_m(\Omega).$ If there exists $\phi\in\mathcal{E}_m(\Omega)$ such that $H_m(\phi)$ is non $m$-polar and if $|u-v|\leq -\phi,$ then $$\ind_{\{u=-\infty\}}H_m(u)=\ind_{\{v=-\infty\}}H_m(v).$$
\end{proposition}

We end this section by recalling a result on the envelope of $m$-subharmonic functions, which plays an important technical role in the proof of the main result of the paper.
\begin{proposition}[Lemma 3.1 in \cite{Pmalay}]\label{promalay}
	Let $\Om$ be a bounded $m$-hyperconvex domain in $\C$ and let $\alpha$ be a finite Radon measure on $\Om$ which puts no mass on $m$-polar sets of $\Om$ such that supp$\alpha\Subset\Om.$ Assume that $\chi: \mathbb{R}^-\to\mathbb{R}^-$ is a nondecreasing continuous function such that $\chi(t)<0$ for all $t<0,$ $\chi(-\infty)=-1$ and $f\in\E_m(\Om)\cap MSH_m(\Om).$ Let $v\in \mathcal F_{m}(f,\Omega)$  be such that $\text{\rm supp}H_{m}(v)\Subset\Omega,$  $H_{m}(v)$ is carried by a $m-$polar set and $\int_{\Om}H_m(v)<+\infty.$ Set 
	$$\mathcal U(\alpha, v)=\left\{\varphi\in SH^{-}_{m}(\Omega): \alpha\leq -\chi (\varphi)H_{m}(\varphi) \text{ and } \varphi\leq v \right\}.$$
	We define
	$$u=(\sup \left\{\varphi: \varphi\in \mathcal U(\alpha, v)\right\})^*.$$
	Then we have $u\in\F_m(f,\Om)$ and
	$-\chi(u)H_m(u)=\alpha+H_m(v).$
\end{proposition}

\section{Proof of the main result}
\n Firstly, we will state and prove the following lemma which is useful in the proof of the main theorem.

\begin{lemma}
	\label{lm1}
	Assume that $g\in\mathcal{E}_m(\Omega)$ and $u\in\mathcal{F}_m(\Omega).$ We consider
	$$h=\sup\{\varphi\in SH^-_m(\Omega):\varphi\leq \min(g,u)\}^*.$$
	Then we have $h\in\mathcal{F}_m(g,\Omega)$ and 
	$$H_m(h)\leq H_m(g)+H_m(u) \,\,\text{on}\,\,\Omega.$$
\end{lemma}

\begin{proof}
	Note that according to Theorem 4.9 in \cite{T19} we have the following assertion: if $u\in\mathcal{F}_m(\Omega)$ then we $\int_{\Omega}H_m(u)<+\infty.$
	Applying Part A of Theorem 3.1 in \cite{PSlova} in the case when $\Omega\equiv\widehat{\Omega}; f\equiv 0$ we get the desired conclusion. 
\end{proof}

\begin{proof}[ The proof of Theorem \ref{main theorem}]
	
	We consider two cases.\\
	\n {\bf Case 1.}  $\chi(-\infty)=-\infty$. Repeating the argument as in Case 2 of Theorem 3.3 in \cite{PDtaiwan}, we get the desired conclusion.\\
	\n {\bf Case 2.}   $\chi(-\infty)>-\infty$.\\
	\n {\sl Step 1}.  In the special case $\chi(-\infty)=-1.$  \\
	We have $$\mu=\ind_{\{w=-\infty\}}\mu+\ind_{\{w>-\infty\}}\mu,$$ where $\mu_1:=\ind_{\{w=-\infty\}}\mu$ is a measure which is carried by an $m$-polar set and $\mu_2:=\ind_{\{w>-\infty\}}\mu$ is a measure which puts no mass on $m$-polar set since $$\ind_{\{w>-\infty\}}\mu\leq \ind_{\{w>-\infty\}}[-\chi(w)]H_m(w).$$
	Note that $$\mu_1\leq \ind_{\{w=-\infty\}}[-\chi(w)]H_m(w)=\ind_{\{w=-\infty\}}[-\chi(-\infty)]H_m(w)\leq H_m(w).$$
	It follows from the hypothesis $\chi(t)<0$ for all $t<0$ and Theorem 3.3 in \cite{DT23} that  $\mathcal{E}_{m,\chi}(\Omega)\subset\mathcal{N}_m(\Omega) \subset\mathcal{E}_m(\Omega).$ Since $w\in\mathcal{E}_{m,\chi}(f,\Omega),$ there exists $\xi\in\mathcal{E}_{m,\chi}(\Omega)\subset\mathcal{E}_m(\Omega)$ such that $f+\xi\leq w\leq f.$ Note that, by Lemma 2.6 in \cite{DT23}, $\mathcal{E}_m(\Omega)$ is a convex cone, we obtain that $f+\xi\in\mathcal{E}_m(\Omega).$ Moreover, $\mathcal{E}_m(\Omega)$ is stable under taking the maximum. Thus, we have $ w=\max(w,f+\xi)\in\mathcal{E}_m(\Omega).$
	According to Theorem 5.9 in \cite{HP17}, there exist a function $v\in\mathcal{E}_m(\Omega) $ such that $$\mu_1=H_m(v) \,\,\text{and}\,\, w\leq v\,\,\text{on}\,\,\Omega.$$
	\n Choosing an increasing sequence $\{\Omega_j\},\Omega_j \Subset \Omega, \Omega_j\nearrow\Omega$ as $j\to\infty.$
	\n Repeating the argument as in Step 3 in the proof of Theorem 5.9 in \cite{HP17}, we can find the decreasing sequence $v_j\in\mathcal{F}_m(\Omega), v_j\geq v$ and $H_m(v_j)=\ind_{\Omega_j}H_m(v).$
	We set
	$$f_j=\sup\{\varphi\in SH_m(\Omega):\varphi\leq \min(f,v_j)\}^*.$$
	By Lemma \ref{lm1} we have $f_j\in\mathcal{F}_m(f,\Omega)$ and \begin{equation}\label{e3.2}H_m(f_j)\leq H_m(f)+H_m(v_j)=H_m(v_j)=\ind_{\Omega_j}H_m(v) \,\,\text{on}\,\,\Omega.\end{equation}
	
	\n Note that $H_m(v)$ is carried by a $m$-polar set. It follows from the inequality \eqref{e3.2} that $H_m(f_j)$ is also carried by a $m$-polar set. Assume that $H_m(f_j)=\ind_{A}H_m(f_j)$ for some $m$-polar set $A.$ Now, we have
	$$H_m(f_j)=\ind_{A}H_m(f_j)= \ind_{A}\ind_{\{f_j=-\infty\}}H_m(f_j)+\ind_{A}\ind_{\{f_j>-\infty\}}H_m(f_j).$$
	By Theorem 2.8 in \cite{DT23}, we see that $\ind_{\{f_j>-\infty\}}H_m(f_j)$ puts no mass on $m$-polar sets. Thus, we have  $\ind_{A}\ind_{\{f_j>-\infty\}}H_m(f_j)$=0.
	This implies that $$H_m(f_j)= \ind_{A}\ind_{\{f_j=-\infty\}}H_m(f_j)=\ind_{\{f_j=-\infty\}}[\ind_{A}H_m(f_j)]=\ind_{\{f_j=-\infty\}}H_m(f_j).$$ 
	Using the same argument, we also have $H_m(v_j)=\ind_{\{v_j=-\infty\}}H_m(v_j).$
	Moreover, we have $$f+v_j\leq f_j\leq v_j.$$  By Proposition \ref{pro2.9}, we infer that
	\begin{equation}\label{e2}H_m(v_j)=\ind_{\{v_j=-\infty\}}H_m(v_j)\leq \ind_{\{f_j=-\infty\}}H_m(f_j)\leq \ind_{\{f+v_j=-\infty\}}H_m(f+v_j).\end{equation}
	Note that $|(f+v_j)-v_j|=-f.$ Moreover, since $f\in MSH_m(\Omega),$ according to Theorem 1.2 in \cite{Bl1}, we have $H_m(f)=0.$ Obviously, $H_m(f)$ puts no mass on $m$-polar sets. By Proposition \ref{pluripolarpart} we obtain that
	\begin{equation}\label{e3} \ind_{\{f+v_j=-\infty\}}H_m(f+v_j)=\ind_{\{v_j=-\infty\}}H_m(v_j)=H_m(v_j)=\ind_{\Omega_j}H_m(v).\end{equation}
	Combining inequality \eqref{e2} with equality \eqref{e3}, we obtain  $H_m(f_j)=\ind_{\{f_j=-\infty\}}H_m(f_j)=\ind_{\Omega_j}H_m(v).$ It means that $H_m(f_j)$ is carried by a $m$-polar set, supp$H_m(f_j)\Subset \Omega.$ Moreover, we also have $f_j\in\mathcal{F}_m(f,\Omega)$ and $\int_{\Omega}H_m(f_j)=\int_{\Omega}\ind_{\Omega_j}H_m(v)\leq \int_{\Omega}H_m(v)< \int_{\Omega}d\mu<\infty.$\\
	Now we will apply  Proposition \ref{promalay}   in the case  $f_j$ plays the role of $v$ and $\alpha_j=\ind_{\{w>-j\}\cap \Omega_j}d\mu_2$ is a finite Radon measure on $\Omega$ which puts no mass on $m$-polar sets such that supp$\alpha_j\Subset\Omega.$\\
	We set 
	$$u_j=(\sup \left\{\varphi: \varphi\in \mathcal U(\alpha_j, f_j)\right\})^*,$$ where
	$$\mathcal U(\alpha_j, f_j)=\left\{\varphi\in SH^{-}_{m}(\Omega): \alpha_j\leq -\chi (\varphi)H_{m}(\varphi) \text{ and } \varphi\leq f_j \right\}.$$
	Then  we have $u_j\in\mathcal{F}_m(f,\Omega)$ and $$-\chi(u_j)H_m(u_j)=\alpha_j+H_m(f_j)\,\,\text{on}\,\,\Omega.$$
	It means we have
	\begin{equation}\label{e3.6}
		-\chi(u_j)H_m(u_j)=\ind_{\{w>-j\}\cap \Omega_j}d\mu_2+\ind_{\Omega_j}H_m(v) \,\,\text{on}\,\,\Omega.
	\end{equation}
	
	\n Since $w\leq v\leq v_j$ and $w\leq f$, we infer that $w\leq f_j$ on $\Omega$. It is easy to see that $w$ belongs to the class of functions which are in the definition of $u_j.$ Thus, we have $w\leq u_j$ on $\Omega$ for all $j.$\\ \n Moreover, since $\{v_j\}_j$ is a decreasing sequence, we also have $\{f_j\}_j$ is a decreasing sequence. Hence, we can check that $\{u_j\}_j$ is also a decreasing sequence. We assume that $u_j\searrow u\geq w\in\mathcal{E}_m(\Omega).$
	According to  Proposition \ref{Pdebrecen}  we have $-\chi(u_j)H_m(u_j)$ converges weakly to $-\chi(u)H_m(u)$ on $\Omega$ as $j\to\infty.$ Moreover, by Lebesgue's dominated convergence theorem, we infer that  $\ind_{\{w>-j\}\cap \Omega_j}d\mu_2 $ converges weakly to $\ind_{\{w>-\infty\}}d\mu_2=d\mu_2$ on $\Omega$ as $j\to\infty.$ By Lebesgue's monotone convergence theorem, we infer that $\ind_{\Omega_j}H_m(v)$ converges weakly to $H_m(v)$  on $\Omega$ as $j\to\infty.$ Thus, let $j\to\infty$ in equality \eqref{e3.6}, we obtain
	\begin{equation*}
		-\chi(u)H_m(u)=d\mu_2+H_m(v)=\mu \,\,\text{on}\,\,\Omega.
	\end{equation*}
	
	Since $w\leq u\leq u_j\leq f_j\leq f$ and $w\in\mathcal{E}_{m,\chi}(f,\Omega)$, we also have $u\in \mathcal{E}_{m,\chi}(f,\Omega).$ \\
	{\sl Step 2}.  In general case $\chi(-\infty)=-M.$ According to the hypothesis $\chi(t)<0$ for all $t<0,$ we infer that $0<M<+\infty.$ Now, we consider $h(t)=\frac{\chi(t)}{M}$. Obviously, we have $h(-\infty)=-1.$ Since $\chi:[-\infty,0]\to [-\infty,0]$ is a nondecreasing continuous function such that $\chi(t)<0$ for all $t<0,$ so is $h(t).$ It follows from the assumption $\mu\leq-\chi(w)H_m(w)$ that $\frac{\mu}{M}\leq -h(w)H_m(w).$ Obviously, we also have $\frac{\mu}{M}$ is a non-negative finite measure on $\Omega.$ Moreover, since $w\in\mathcal{E}_{m,\chi}(f,\Omega),$ there exists a fucntion $\xi\in\mathcal{E}_{m,\chi}(\Omega)$ such that $f+\xi\leq w\leq f.$ By the definition of the class $\mathcal{E}_{m,\chi}(\Omega)$, there exists a decreasing function $(\xi_j)\in\mathcal{E}_m^0(\Omega),\xi_j\searrow\xi$ such that $\sup\limits_{j}\int\limits_{\Omega}-\chi(\xi_j)H_m(\xi)<+\infty.$ Thus, we also have  $\sup\limits_{j}\int\limits_{\Omega}-\frac{\chi(\xi_j)}{M}H_m(\xi)<+\infty.$ This means we have $\sup\limits_{j}\int\limits_{\Omega}-h(\xi_j)H_m(\xi)<+\infty.$ Therefore, we obtain that $\xi\in \mathcal{E}_{m,h}(\Omega).$ Hence, we infer that $w\in\mathcal{E}_{m,h}(f,\Omega).$ Applying Step 1 with $\frac{\mu}{M}$ plays the role of $\mu$ and $h$ plays the role of $\chi,$ there exists a function $u\in \mathcal{E}_{m,h}(f,\Omega)$ such that $u\geq w$ and $-h(u)H_m(u)=\frac{\mu}{M}.$ This implys that $-\chi(u)H_m(u)=\mu.$ It is easy to see that $u\in \mathcal{E}_{m,\chi}(f,\Omega).$
	The proof is complete.
\end{proof}

\section*{Declarations}
\subsection*{Ethical Approval}
This declaration is not applicable.
\subsection*{Competing interests}
The author has no conflicts of interest to declare that are relevant to the content of this article.

\subsection*{Availability of data and materials}
This declaration is not applicable.

\end{document}